\documentclass[11pt]{amsart}
\usepackage{stmaryrd}
\expandafter\def\csname opt@stmaryrd.sty\endcsname
{only,shortleftarrow,shortrightarrow}
\usepackage{extpfeil}
\usepackage{amsmath,amssymb,amsthm}
\usepackage{aliascnt}
\usepackage{hyperref}
\hypersetup{colorlinks=true,urlcolor=blue,citecolor=blue,linkcolor=blue}
\usepackage{courier}
\usepackage{tikz}
\usepackage{tikz-cd}
\usetikzlibrary{calc,matrix,arrows,decorations.markings}
\usepackage{array}
\usepackage{color}
\usepackage{enumerate}
\usepackage{nicefrac}
\usepackage{listings}
\usepackage{siunitx}
\usepackage{seqsplit}
\usepackage{longtable}
\usepackage{colonequals}
\newcommand{\Q}{\mathbb{Q}}

\newcommand{\R}{\mathbb{R}}
\newcommand{\C}{\mathbb{C}}
\renewcommand{\H}{\mathbb{H}}
\newcommand{\Z}{\mathbb{Z}}

\renewcommand{\P}{\mathbb{P}}

\newcommand{\End}{\operatorname{End}}

\newcommand{\gon}{\operatorname{gon}}

\newcommand{\CM}{\operatorname{CM}}

\DeclareMathOperator{\PSL}{PSL}

\DeclareFontFamily{U}{wncy}{}
    \DeclareFontShape{U}{wncy}{m}{n}{<->wncyr10}{}
    \DeclareSymbolFont{mcy}{U}{wncy}{m}{n}
    \DeclareMathSymbol{\Sh}{\mathord}{mcy}{"58}

\newtheorem{theorem}{Theorem}[section]
\theoremstyle{definition}

\newtheorem{remark}{Remark}[section]
\newtheorem{proposition}{Proposition}[section]

\newtheorem{lemma}{Lemma}[section]

\usepackage{xcolor}

\begin{document}

\title{Shimura curves admitting a smooth plane model}

\author{\sc Oana Padurariu}
\address{Oana Padurariu \\
Max-Planck-Institut für Mathematik Bonn\\
Germany}
\urladdr{https://sites.google.com/view/oanapadurariu/home}
\email{oana.padurariu11@gmail.com}

\maketitle
\begin{abstract}
 We prove that there are no Shimura curves $X_0^D(N)$ with $N$ squarefree of genus $g \ge 2$  admitting a smooth plane model.
\end{abstract}


\section{Introduction}
Modular curves and Shimura curves are central moduli spaces in number theory. 
Modular curves parametrize elliptic curves equipped with level structure and Shimura curves parametrize abelian surfaces with quaternionic multiplication. Shimura curves are naturally compact, in contrast to modular curves that are compactified by adding cusps. These cusps are extremely useful in the computations of defining equations for modular curves, which lead to an impressive database of such explicit equations.
On the other hand, there is a scarcity of defining equations for Shimura curves in the literature.

We briefly recall some results
about the classification of Shimura curves. Gonz\'{a}lez--Rotger \cite{GR06} provided defining equations for non-elliptic Shimura curves $X_0^D(N)$ of genus one. Voight \cite{Voight09} fully classified Shimura curves of genus at most two. Guo--Yang \cite{GY17} gave a complete list of geometrically hyperelliptic Shimura curves $X_0^D(N)$ with explicit equations. The next step in the classification of Shimura curves is finding equations for all Shimura curves admitting a smooth plane model, which is what Anni--Assaf--Lorenzo-Garc\'{\i}a \cite{AALG23} achieved for modular curves. A natural question, inspired by \cite{AALG23}, is whether there are any Shimura curves $X_0^D(N)$ admitting a smooth plane model beyond those of genus at most one. In this paper, we prove the following:

\begin{theorem}\label{theorem:main_theorem}
Let $X_0^D(N)$ be a Shimura curve for which the level $N$ is squarefree. If $X_0^D(N)$ admits a smooth plane model, then the genus of $X_0^D(N)$ is at most one.
\end{theorem}

The proof relies on the interplay between the genus and the gonality of Shimura curves, and on computations of the number of CM points fixed by Atkin--Lehner involutions.  
The paper is accompanied by code \cite{MyCode} written in Magma \cite{magma}. 

\subsection*{Acknowledgements}
It is a pleasure to thank Pete L. Clark, Ciaran Schembri, Frederick Saia, and John Voight for many helpful conversations about Shimura curves. 
The author is very grateful to Max-Planck-Institut für Mathematik Bonn for hospitality and financial support. 

\section{Some background on Shimura curves}
Throughout the paper, $p\in\Z_{\ge 1}$ is a prime number. In this section, we recall some definitions and results about Shimura curves.

\subsection{Quaternion algebras}
To define Shimura curves, one first needs to define quaternion algebras, which are extensively covered in \cite{VoightQA}. A \textit{quaternion algebra} $B$ defined over $\Q$ is a central simple $\Q$-algebra of dimension $4$. Alternatively, $B$ is a quaternion algebra defined over $\Q$ if and only if there are elements $i,j \in B$ such that
\[
B \simeq \Q \cdot 1 + \Q \cdot i + \Q \cdot j + \Q \cdot ij 
\]
where $i^2,j^2 \in \Q^{\times}$ with $ij = -ji$. An example of a quaternion algebra is $M_2(\Q)$.
An \textit{order} of $B$ is a full $\Z$-lattice which is also a subring.

A very important invariant for a quaternion algebra is the \textit{discriminant}. For any place $v$ of $\Q$, we say that $B$ is \textit{ramified} at $v$ if the completion $B \otimes_{\Q} \Q_v$ is a division algebra. The set of ramified places of $B$ is finite, of even cardinality, and uniquely determines the isomorphism class of a quaternion algebra. Then the discriminant $D$ of $B$ is defined as
\[
  D \colonequals \prod_{\substack{v\text{ finite} \\ \text{ramified}}} v.  
\]

We say that $B$ is \textit{indefinite} if $B$ is unramified at the archimedean place. Note that in this case $\omega(D)$ is even, where 
\[
\omega(n) \colonequals \#\{ p : p |n  \}.
\]

We say that $d$ is a \textit{Hall divisor} of $n$ if $d|n$ and $\gcd(d,n/d) = 1$. We denote this by $d \parallel n$.
\smallskip

\subsection{Definition of Shimura curves}
Let $B$ be an indefinite quaternion algebra defined over $\Q$ of discriminant $D$ and let $O$ denote a maximal order. There are embeddings:
\begin{align*}
\iota_\infty &: B \hookrightarrow M_2(\R) \simeq B \otimes_{\Q} \R \\
\iota_p &: O \hookrightarrow M_2(\Z_p) \simeq O \otimes_{\Z} \Z_p
\end{align*}
for any prime $p$ not dividing $D$. We have a notion of level, similar to the case of modular curves: for any $N \in \Z_{\ge 1}$ coprime to $D$, we define the \textit{Eichler order} of level $N$ to be 
\[
O_N  \colonequals \{ x \in O \ | \ \iota_p(x) \equiv \begin{pmatrix} \ast & \ast  \\ 0 & \ast \end{pmatrix} \ \text{mod} \ p^e \ \text{for all} \ p^e \parallel N \}.
\]

Let $O_N^1$ be the set of norm one elements of $O_N$. Then $O_N^1$ acts on the upper half plane $\H = \{ x + iy \ | \ y > 0 \}$ via $\iota_{\infty}(O_N^1)/\{ \pm 1 \} \le \PSL_2(\R)$ by 

$$ \gamma \cdot z \colonequals \iota_\infty(\gamma) \cdot z = \begin{pmatrix} a & b \\ c& d \end{pmatrix} \cdot z = \frac{az+b}{cz+d}.$$

Consider the Riemann surface
$$\left(\iota_\infty (O_N^1) / \{ \pm 1 \} \right)\backslash \H.$$
If $B$ has discriminant $D=1$, then this Riemann surface is the classical modular curve $Y_0(N).$

From now on we shall assume that $B$ is not isomorphic to $M_2(\Q)$, or equivalently that $D>1.$  By work of Shimura \cite[Main Theorem I]{Shimura67}, \cite[Theorem 2.5]{Shimura70}, there is an algebraic curve  \[X_0^D(N)\] such that there is an open immersion of Riemann surfaces \[\left(\iota_\infty (O_N^1) / \{ \pm 1 \} \right)\backslash \H \hookrightarrow X_0^D(N)\] which is a biregular isomorphism.  The curve $X_0^D(N)$ has a canonical model defined over $\Q$ \cite[\S3]{Shimura67} and it is called a \textit{Shimura curve} of discriminant $D$ and level $N$.  

As with modular curves,  Shimura curves arise naturally as a moduli problem.  The curve $X_0^D(N)$ parameterizes pairs $(A,\iota)$, where 
\begin{itemize}
\item $A$ is an abelian surface, and
\item $\iota : O_N \hookrightarrow \End(A)$ is an embedding.
\end{itemize}

We call such a pair $(A,\iota)$ a \textit{QM-surface}.  It is well known that Shimura curves have no real points \cite[Theorem 0]{Shimura75}, \cite[\S3]{Ogg85},  i.e.
$$X_0^D(1)(\R) = \emptyset.$$

\subsection{Atkin--Lehner involutions}
As in the case of modular curves, Shimura curves also have Atkin--Lehner involutions.
Consider the normalizer group
$$N_{B_{>0}^\times}(O_N^1) \colonequals \{  \alpha \in B_{>0}^\times \ | \ \alpha^{-1} O_N^1 \alpha = O_N^1  \},$$
where $B_{>0}^\times$ are the units in $B$ of positive norm.  Elements of $N_{B_{>0}^\times}(O_N^1)$ naturally define an automorphism of the Shimura curve $X_0^D(N)$,  with the scalars $\Q^\times$ acting trivially.  It is straightforward to check that two elements $\alpha, \beta$ of $N_{B_{>0}^\times}(O_N^1)$ define the same automorphism if and only if $\alpha \beta^{-1} \in \Q^\times O_N^1.$ We define the group of Atkin--Lehner involutions as follows:
$$W(D,N) \colonequals N_{B_{>0}^\times}(O_N^1) / \Q^\times O_N^1.$$
We have an Atkin--Lehner involution $w_m$ for every Hall divisor $m$ of $DN$.
There is an identification \[W(D,N) \simeq (\Z/2\Z)^{\omega(DN)}\] and we write $W(D, N) = \{ w_m \ | \ m \parallel ND \}$ \cite[Ch.  28]{VoightQA}. 

\subsection{Genus of Shimura curves}

In this subsection we recall an explicit formula for the genus of $X_0^D(N)$.
For coprime positive integers $D$ and $N$ such that $D$ is squarefree, define

\[
e_k(D,N) \colonequals  \prod_{p|D} \left(1 - \left(\frac{-k}{p}\right)\right) \prod_{q \parallel N} \left(1 + \left(\frac{-k}{q}\right)\right) \prod_{q^2|N}\nu_p(k) ,
\]
where
\[   
\nu_p(k) = 
     \begin{cases}
       2 &\quad\text{if } \left(\frac{-k}{q}\right) = 1, \\
       0 &\quad\text{otherwise},  \\
     \end{cases}
\]
and $\left(\frac{\cdot}{\cdot}\right)$
is the Kronecker quadratic symbol.

\begin{remark}
The above quantities have geometric meaning, for example:
    \begin{itemize}
        \item $e_4(D,N)$ is the number of $\Z[i]$-CM points on $X_0^D(N)$, and also the number of order $2$ elliptic points on $X_0^D(N)$,
        \item $e_3(D,N)$ is the number of $\Z\left[\frac{1+\sqrt{-3}}{2}\right]$-CM points on $X_0^D(N)$, and also the number of order $3$ elliptic points on $X_0^D(N)$.
    \end{itemize}
\end{remark}

Let $\varphi$ and $\psi$ be the multiplicative functions such that 

\[
\varphi(p^k) = p^k - p^{k-1}, \quad \psi(p^k) = p^k + p^{k-1}.
\]

\begin{proposition}\cite[p.280, 301]{Ogg83}
For $D > 1$, the genus of $X_0^D(N)$ is
\[
g(X_0^D(N)) = 1 + \frac{\varphi(D)\psi(N)}{12} - \frac{e_4(D,N)}{4} - \frac{e_3(D,N)}{3}.
\]
\end{proposition}

\begin{remark}
A Shimura curve $X_0^D(N)$ with $N$ squarefree is of genus zero if and only if the pair $(D,N)$ is in the set:
\[
\{ (6,1),(10,1),(22,1) \}.
\]
A Shimura curves $X_0^D(N)$ with $N$ squarefree is of genus one if and only if the pair $(D,N)$ is in the set:
\[
\{ (14, 1),(15, 1),(21, 1),(6, 5),(10, 3),(33, 1),(34, 1),(6, 7),(46, 1),(10, 7),(6, 13)\}.
\]
Equations for these genus one curves and formulae for their Atkin--Lehner involutions can be found in \cite[Theorem 3.4]{GR06}.
\end{remark}

A lower bound for the genus is given by the following: 
\begin{lemma}{\cite[Lemma 10.5]{Saia22}}
\label{lemma: genus_upper_bound}
For $D>1$ an indefinite rational quaternion discriminant and $N \in \mathbb{Z}_{\ge 1}$ relatively prime to $D$, we have
\[ g(X_0^D(N)) > 1 + \frac{DN}{12}\left( \frac{1}{e^\gamma \log\log(DN) + \frac{3}{\log\log{6}}} \right)- \frac{7\sqrt{DN}}{3},\]
where $\gamma$ is the Euler constant.
\end{lemma}

\subsection{Gonality of a curve}
For a field $K$ and a curve $C/K$, the $K$-\textit{gonality} $\gon_K(C)$ is defined to be the least degree of a non-constant morphism $f : C \rightarrow \P^1$ defined over $K$. For $\overline{K}$ an algebraic closure of $K$, the $\overline{K}$-gonality is also called the \textit{geometric gonality}. 
We will use the following result to obtain an upper bound on the genus of $X_0^D(N)$:

\begin{theorem}{\cite[Theorem 1.1]{Abramovich96}} \label{theorem: genus_gonality_inequality}
\[
\frac{21}{200} \left( g(X_0^D(N))-1\right) \le \gon_{\C} (X_0^D(N)).
\]   
\end{theorem}

\subsection{CM points}
There is a special type of points on Shimura curve called \textit{CM points} that play a very important rule. As Shimura curves lack cusps, Guo--Yang \cite{GY17} used CM points to help in the computation of defining equations for all geometrically hyperelliptic Shimura curves $X_0^D(N)$. Moreover, if a point is fixed by some Atkin--Lehner involution, then it is a CM point.
This observation is one of the main ingredients in the proof of \autoref{theorem:main_theorem}.

Let $K = \Q(\sqrt{-d})$ be an imaginary quadratic field which admits an embedding \[q : K \hookrightarrow B.\]  Let $R \subset K$ be an order.  Then $R$ is said to be \textit{optimally embedded} in $O_N$ if we have $q(K) \cap O_N = q(R).$

For an optimally embedded order $R$, there is exactly one fixed point $z$ of the quotient $\left(\iota_\infty (O_N^1) / \{ \pm 1 \} \right)\backslash \H$ under the action of $\iota_{\infty}( q(R)).$ We say that $z$ is a \textit{CM point} for $R$.  The set 
$$\CM(R) := \{ z \in \left(\iota_\infty (O_N^1) / \{ \pm 1 \} \right)\backslash \H \ | \ \iota_{\infty}( q(R)) \cdot z = z, \ q : R \hookrightarrow B \ \text{an optimal embedding} \}$$
is the set of CM points for $R$. We identify the CM points on $\left(\iota_\infty (O_N^1) / \{ \pm 1 \} \right)\backslash \H$ with the corresponding points on $X_0^D(N)$.  

For any $W \le W(D,N)$, there is the natural quotient map \[\pi_W : X_0^D(N) \longrightarrow X_0^D(N)/W.\] We say a point $Q \in X_0^D(N)/W$ is a \textit{CM point} of $X_0^D(N)/W$ if any of its preimages on $X_0^D(N)$ is a CM point. 

\begin{remark} 
CM points on $X_0^D(N)$ correspond to abelian surfaces with extra endomorphisms by the CM field $K$: \[\End(A) \otimes_\Z \Q \simeq M_2(K).\] 

We state an alternative description of CM points.
Recall the cover map \[X_0^D(N)\to X_0^D(1).\]
A CM point for $R$ on $X_0^D(N)$ is the same as a point in the fibre over $[(A,\iota)] \in X_0^D(1)$ such that \[\End(A) \otimes_\Z \Q \simeq M_2(K)\] and the ring of QM-equivariant endomorphisms of $A$ is isomorphic to $R$ (cf. \cite[\S 2.2]{Saia22}).
\end{remark}

\begin{theorem}
\cite[Section 1]{Ogg83}
\label{theorem: CM_fixed_by_AL}
Let $N$ be squarefree and $m|DN$.   
The set of fixed points of $X_0^D(N)$ under the involution $w_m$ is
$$X_0^D(N)^{\langle w_m \rangle} =
\begin{cases}
\CM(\Z[\sqrt{-1}]) \cup \CM(\Z[\sqrt{-2}]) &\text{if} \ m=2, \\
\CM(\Z[\sqrt{-m}]) \cup \CM\left(\Z\left[\frac{1+\sqrt{-m}}{2}\right]\right) &\text{if} \ m \equiv 3 \ \text{mod} \ 4, \\
\CM(\Z[\sqrt{-m}]) \ &\text{otherwise}.
\end{cases}
 $$
\end{theorem}

We would like to compute $\#X_0^D(N)^{\langle w_m \rangle}$. To this end, we need to first understand $\#CM(R)$ for an order $R$.

Assume $N$ is squarefree. We will use the following notation:

\begin{itemize}
    \item $\Delta_R$ be the discriminant of $R$,
    \item $f$ the conductor of $R$,
    \item $h(R)$ the class number of $R$,
    \item $\big( \frac{K}{p} \big)$ be the Kronecker symbol,
    \item  $\big( \frac{R}{p} \big) \colonequals \big( \frac{K}{p} \big)$ if $p \nmid f$ and $1$ otherwise.
\end{itemize}

The following three quantities are needed in the computation of $\#\CM(R)$:
$$D(R) \colonequals \prod_{p |D, \ \big( \frac{R}{p} \big) = -1} p, \ \ \ N(R) \colonequals \prod_{p|N, \ \big( \frac{R}{p} \big) = 1} p, \ \ \ N^\ast(R) \colonequals \prod_{p|N,  \ p \nmid f, \big( \frac{R}{p} \big)=1} p.$$

\begin{theorem}[{\cite[Section 1]{Ogg83}}, \cite{Jordan81}, {\cite[Lemma 2.5]{BD96}}]
\label{theorem: number_of_CM_points}
The set $\CM(R)$ is nonempty if and only if $\frac{DN}{D(R)N^\ast(R)} | \Delta_R$ and in this case 
\[
\# \CM(R)  = 2^{\omega(D(R) N(R))} \cdot h(R).
\]

\end{theorem}

The code that computes $\# \CM(R)$ and $\# X_0^D(N)^{\langle w_m \rangle}$ was previously implemented as part of the code accompanying \cite{PS23}.

\section{Proof of \autoref{theorem:main_theorem}}

Assume $X_0^D(N)$ is a Shimura curve with squarefree level $N$ that admits a degree $d$ smooth plane model over $\Q$. We know that the geometric gonality of a smooth plane curve of degree $d$ is $d-1$, see \cite[Theorem A]{CK90}, while the genus-degree formula states that $g(X_0^D(N)) = \frac{(d-1)(d-2)}{2}$. Then using \autoref{theorem: genus_gonality_inequality} we obtain:
\[
\frac{21}{200} \left( \frac{(d-1)(d-2)}{2}
-1\right) = \frac{21}{200} ( g(X_0^D(N))-1) \le \gon_{\C} (X_0^D(N)) = d-1,
\]
which implies that $d \le 21$. Thus the genus of $X_0^D(N)$ is at most $190$.

By Lemma \ref{lemma: genus_upper_bound}, we find that for $DN > 110011$, we have $g(X_0^D(N)) > 190$. Therefore, we only need to consider the pairs $(D,N)$ where $DN \le 110011$ is squarefree. We may consider only those pairs for which $\omega(D)$ is even, as the set of ramified places has even cardinality. We further need that 
\[
g(X_0^D(N)) \in \left \{ \frac{(d-1)(d-2)}{2} : 1 \le d \le 21  \right \}.
\]
Thus we are left to analyze $312$ pairs $(D,N)$. The discriminant $D$ is at most $6990$, while the level $N$ is at most $1033$.
Let $(D,N)$ be such a pair, for which \[g(X_0^D(N)) = \frac{(d-1)(d-2)}{2}\] for some $d \in \{ 1, \cdots, 21 \}$.

By \autoref{theorem: CM_fixed_by_AL}, only CM points can be fixed by an Atkin--Lehner involution $w_m$. Further, by \autoref{theorem: number_of_CM_points}, one is able to compute $\#X_0^D(N)^{\langle w_m \rangle}$. We will combine these results with the following:

\begin{theorem}{\cite[Remark 2.1 (i) and Theorem 2.2 with $n=2$]{HKKO10}} \label{theorem: fixed_points}
Let $C$ be a smooth plane curve of degree $d \ge 4$ and $\sigma$ an involution of $C$. Then the involution $\sigma$ has $f = d + \frac{1-(-1)^d}{2}$ fixed points.
\end{theorem}

If $X_0^D(N)$ admits a plane model of degree $d \ge 4$, then by \autoref{theorem: fixed_points} we have that:
\begin{equation}\label{eqn: equal_genera}
\#X_0^D(N)^{\langle w_m \rangle} = d + \frac{1-(-1)^d}{2} \quad \text{for all } \, 1 < m \parallel DN.
\end{equation}
For $298$ out of the $312$ candidate pairs $(D,N)$ we have that $g(X_0^D(N)) \ge 3$, while for the remaining $14$ we have $g(X_0^D(N)) \le 1$.
Running our code \cite{MyCode} on these $298$ pairs $(D,N)$ we find that Property (\ref{eqn: equal_genera}) does not hold for any of the curves with $g(X_0^D(N)) \ge 3$. Thus there is no Shimura curve $X_0^D(N)$ with $N$ squarefree admitting a smooth plane model of degree $d \ge 4$.
In conclusion, if $X_0^D(N)$ admits a smooth plane model for some $N$ squarefree, then $g(X_0^D(N)) \le 1$.
\setcounter{tocdepth}{1}

\bibliographystyle{amsalpha}
\bibliography{biblio}

\providecommand{\bysame}{\leavevmode\hbox to3em{\hrulefill}\thinspace}
\providecommand{\MR}{\relax\ifhmode\unskip\space\fi MR }
\providecommand{\MRhref}[2]{%
  \href{http://www.ams.org/mathscinet-getitem?mr=#1}{#2}
}
\providecommand{\href}[2]{#2}
\begin{thebibliography}{HKKO10}

\bibitem[AALG23]{AALG23}
Samuele Anni, Eran Assaf, and Elisa Lorenzo~Garc\'{\i}a, \emph{On smooth plane models for modular curves of {S}himura type}, Res. Number Theory \textbf{9} (2023), no.~2, Paper No. 21, 20. \MR{4563688}

\bibitem[Abr96]{Abramovich96}
Dan Abramovich, \emph{A linear lower bound on the gonality of modular curves}, Internat. Math. Res. Notices (1996), no.~20, 1005--1011. \MR{1422373}

\bibitem[BCP97]{magma}
W.~Bosma, J.~Cannon, and C.~Playoust, \emph{The {M}agma algebra system. {I}. {T}he user language ({M}agma {V}2.26-5)}, J. Symbolic Comput. \textbf{24} (1997), no.~3--4, 235--265, Computational algebra and number theory (London, 1993).

\bibitem[BD96]{BD96}
M.~Bertolini and H.~Darmon, \emph{Heegner points on {M}umford-{T}ate curves}, Invent. Math. \textbf{126} (1996), no.~3, 413--456. \MR{1419003}

\bibitem[CK90]{CK90}
Marc Coppens and Takao Kato, \emph{The gonality of smooth curves with plane models}, Manuscripta Math. \textbf{70} (1990), no.~1, 5--25. \MR{1080899}

\bibitem[GR06]{GR06}
Josep Gonz\'{a}lez and Victor Rotger, \emph{Non-elliptic {S}himura curves of genus one}, J. Math. Soc. Japan \textbf{58} (2006), no.~4, 927--948. \MR{2276174}

\bibitem[GY17]{GY17}
Jia-Wei Guo and Yifan Yang, \emph{Equations of hyperelliptic {S}himura curves}, Compos. Math. \textbf{153} (2017), no.~1, 1--40. \MR{3622871}

\bibitem[HKKO10]{HKKO10}
Takeshi Harui, Takao Kato, Jiryo Komeda, and Akira Ohbuchi, \emph{Quotient curves of smooth plane curves with automorphisms}, Kodai Math. J. \textbf{33} (2010), no.~1, 164--172. \MR{2732237}

\bibitem[Jor81]{Jordan81}
Bruce~Winchester Jordan, \emph{On the diophantine arithmetic of {S}himura curves}, ProQuest LLC, Ann Arbor, MI, 1981, Thesis (Ph.D.)--Harvard University. \MR{2936886}

\bibitem[Ogg83]{Ogg83}
A.~P. Ogg, \emph{Real points on {S}himura curves}, Arithmetic and geometry, {V}ol. {I}, Progr. Math., vol.~35, Birkh\"{a}user Boston, Boston, MA, 1983, pp.~277--307. \MR{717598}

\bibitem[Ogg85]{Ogg85}
\bysame, \emph{Mauvaise r\'{e}duction des courbes de {S}himura}, S\'{e}minaire de th\'{e}orie des nombres, {P}aris 1983--84, Progr. Math., vol.~59, Birkh\"{a}user Boston, Boston, MA, 1985, pp.~199--217. \MR{902833}

\bibitem[Pad]{MyCode}
Oana Padurariu, \emph{{M}agma code}, \url{https://github.com/oana-adascalitei/SmoothPlane}.

\bibitem[PS23]{PS23}
Oana Padurariu and Ciaran Schembri, \emph{Rational points on {A}tkin--{L}ehner quotients of geometrically hyperelliptic {S}himura curves}, Expo. Math. \textbf{41} (2023), no.~3, 492--513.

\bibitem[Sai22]{Saia22}
F.~Saia, \emph{{CM points on Shimura Curves via QM-equivariant isogeny volcanoes}}, arXiv:2212.12635 (2022).

\bibitem[Shi67]{Shimura67}
Goro Shimura, \emph{Construction of class fields and zeta functions of algebraic curves}, Ann. of Math. (2) \textbf{85} (1967), 58--159. \MR{0204426}

\bibitem[Shi70]{Shimura70}
\bysame, \emph{On canonical models of arithmetic quotients of bounded symmetric domains}, Ann. of Math. (2) \textbf{91} (1970), 144--222. \MR{257031}

\bibitem[Shi75]{Shimura75}
\bysame, \emph{On the real points of an arithmetic quotient of a bounded symmetric domain}, Math. Ann. \textbf{215} (1975), 135--164. \MR{0572971}

\bibitem[Voi09]{Voight09}
John Voight, \emph{Shimura curves of genus at most two}, Math. Comp. \textbf{78} (2009), no.~266, 1155--1172. \MR{2476577}

\bibitem[Voi21]{VoightQA}
\bysame, \emph{Quaternion algebras}, Graduate Texts in Mathematics, vol. 288, Springer, Cham, [2021] \copyright 2021. \MR{4279905}

\end{thebibliography}
\end{document}